\documentclass[12pt]{article}

\usepackage{amsfonts,amssymb,amsthm,amsmath,latexsym}
\usepackage{cite}

\begin{document}

\title{A really simple elementary proof \\
       of the uniform boundedness theorem}

\author{
     {\small Alan D.~Sokal\thanks{Also at Department of Mathematics,
           University College London, London WC1E 6BT, England.}}  \\[-2mm]
     {\small\it Department of Physics}       \\[-2mm]
     {\small\it New York University}         \\[-2mm]
     {\small\it 4 Washington Place}          \\[-2mm]
     {\small\it New York, NY 10003 USA}      \\[-2mm]
     {\small\tt sokal@nyu.edu}               \\[-2mm]
     {\protect\makebox[5in]{\quad}}  
     \\
}

\date{May 7, 2010 \\[3mm]
      revised October 20, 2010 \\[1mm]
      to appear in the {\em American Mathematical Monthly}}

\maketitle
\thispagestyle{empty}   
\begin{abstract}
I give a proof of the uniform boundedness theorem that is elementary
(i.e., does not use any version of the Baire category theorem)
and also extremely simple.
\end{abstract}


\bigskip
\noindent
{\bf Mathematics Subject Classification (MSC 2000) codes:}
46B99 (Primary); 46B20, 46B28 (Secondary).

\clearpage

\newtheorem{theorem}{Theorem}[section]
\newtheorem{proposition}[theorem]{Proposition}
\newtheorem{lemma}[theorem]{Lemma}
\newtheorem{corollary}[theorem]{Corollary}
\newtheorem{definition}[theorem]{Definition}
\newtheorem{conjecture}[theorem]{Conjecture}
\newtheorem{question}[theorem]{Question}
\newtheorem{example}[theorem]{Example}
\newtheorem{remark}[theorem]{Remark}

\renewcommand{\theenumi}{\alph{enumi}}
\renewcommand{\labelenumi}{(\theenumi)}
\def\eop{\hbox{\kern1pt\vrule height6pt width4pt
depth1pt\kern1pt}\medskip}
\def\prf{\par\noindent{\bf Proof.\enspace}\rm}
\def\rmk{\par\medskip\noindent{\bf Remark\enspace}\rm}

\newcommand{\be}{\begin{equation}}
\newcommand{\ee}{\end{equation}}
\newcommand{\<}{\langle}
\renewcommand{\>}{\rangle}
\newcommand{\widebar}{\overline}
\def\reff#1{(\protect\ref{#1})}
\def\spose#1{\hbox to 0pt{#1\hss}}
\def\ltapprox{\mathrel{\spose{\lower 3pt\hbox{$\mathchar"218$}}
    \raise 2.0pt\hbox{$\mathchar"13C$}}}
\def\gtapprox{\mathrel{\spose{\lower 3pt\hbox{$\mathchar"218$}}
    \raise 2.0pt\hbox{$\mathchar"13E$}}}
\def\textprime{${}^\prime$}
\def\proof{\par\medskip\noindent{\sc Proof.\ }}
\def\firstproof{\par\medskip\noindent{\sc First Proof.\ }}
\def\secondproof{\par\medskip\noindent{\sc Second Proof.\ }}
\def\thirdproof{\par\medskip\noindent{\sc Third Proof.\ }}
\def\qed{ $\square$ \bigskip}
\def\proofof#1{\bigskip\noindent{\sc Proof of #1.\ }}
\def\firstproofof#1{\bigskip\noindent{\sc First Proof of #1.\ }}
\def\secondproofof#1{\bigskip\noindent{\sc Second Proof of #1.\ }}
\def\thirdproofof#1{\bigskip\noindent{\sc Third Proof of #1.\ }}
\def\half{ {1 \over 2} }
\def\third{ {1 \over 3} }
\def\twothird{ {2 \over 3} }
\def\smfrac#1#2{{\textstyle{#1\over #2}}}
\def\smhalf{ {\smfrac{1}{2}} }
\newcommand{\real}{\mathop{\rm Re}\nolimits}
\renewcommand{\Re}{\mathop{\rm Re}\nolimits}
\newcommand{\imag}{\mathop{\rm Im}\nolimits}
\renewcommand{\Im}{\mathop{\rm Im}\nolimits}
\newcommand{\sgn}{\mathop{\rm sgn}\nolimits}
\newcommand{\tr}{\mathop{\rm tr}\nolimits}
\newcommand{\supp}{\mathop{\rm supp}\nolimits}
\def\hboxscript#1{ {\hbox{\scriptsize\em #1}} }
\renewcommand{\emptyset}{\varnothing}

\newcommand{\restrict}{\upharpoonright}

\newcommand{\scra}{{\mathcal{A}}}
\newcommand{\scrb}{{\mathcal{B}}}
\newcommand{\scrc}{{\mathcal{C}}}
\newcommand{\scre}{{\mathcal{E}}}
\newcommand{\scrf}{{\mathcal{F}}}
\newcommand{\scrg}{{\mathcal{G}}}
\newcommand{\scrh}{{\mathcal{H}}}
\newcommand{\scrk}{{\mathcal{K}}}
\newcommand{\scrl}{{\mathcal{L}}}
\newcommand{\scro}{{\mathcal{O}}}
\newcommand{\scrp}{{\mathcal{P}}}
\newcommand{\scrr}{{\mathcal{R}}}
\newcommand{\scrs}{{\mathcal{S}}}
\newcommand{\scrt}{{\mathcal{T}}}
\newcommand{\scrv}{{\mathcal{V}}}
\newcommand{\scrw}{{\mathcal{W}}}
\newcommand{\scrz}{{\mathcal{Z}}}

\newcommand{\ahat}{{\widehat{a}}}
\newcommand{\Zhat}{{\widehat{Z}}}
\renewcommand{\k}{{\mathbf{k}}}
\newcommand{\n}{{\mathbf{n}}}
\newcommand{\vv}{{\mathbf{v}}}
\newcommand{\bv}{{\mathbf{v}}}
\newcommand{\w}{{\mathbf{w}}}
\newcommand{\x}{{\mathbf{x}}}
\newcommand{\g}{{\boldsymbol{g}}}
\newcommand{\cc}{{\mathbf{c}}}
\newcommand{\zero}{{\mathbf{0}}}
\newcommand{\one}{{\mathbf{1}}}
\newcommand{\balpha}{{\boldsymbol{\alpha}}}

\newcommand{\C}{{\mathbb C}}
\newcommand{\D}{{\mathbb D}}
\newcommand{\Z}{{\mathbb Z}}
\newcommand{\N}{{\mathbb N}}
\newcommand{\Q}{{\mathbb Q}}
\newcommand{\R}{{\mathbb R}}
\newcommand{\RR}{{\mathbb R}}

\def\cbar{{\overline{C}}}
\def\ctilde{{\widetilde{C}}}
\def\zbar{{\overline{Z}}}
\def\pitilde{{\widetilde{\pi}}}

%
%
\newcommand{\stirlingsubset}[2]{\genfrac{\{}{\}}{0pt}{}{#1}{#2}}
\newcommand{\stirlingcycle}[2]{\genfrac{[}{]}{0pt}{}{#1}{#2}}
\newcommand{\assocstirlingsubset}[3]{{\genfrac{\{}{\}}{0pt}{}{#1}{#2}}_{\! \ge #3}}
\newcommand{\assocstirlingcycle}[3]{{\genfrac{[}{]}{0pt}{}{#1}{#2}}_{\ge #3}}
\newcommand{\genstirlingsubset}[4]{{\genfrac{\{}{\}}{0pt}{}{#1}{#2}}_{\! #3,#4}}
\newcommand{\euler}[2]{\genfrac{\langle}{\rangle}{0pt}{}{#1}{#2}}
\newcommand{\eulergen}[3]{{\genfrac{\langle}{\rangle}{0pt}{}{#1}{#2}}_{\! #3}}
\newcommand{\eulersecond}[2]{\left\langle\!\! \euler{#1}{#2} \!\!\right\rangle}
\newcommand{\eulersecondgen}[3]{{\left\langle\!\! \euler{#1}{#2} \!\!\right\rangle}_{\! #3}}
\newcommand{\binomvert}[2]{\genfrac{\vert}{\vert}{0pt}{}{#1}{#2}}


\newenvironment{sarray}{
             \textfont0=\scriptfont0
             \scriptfont0=\scriptscriptfont0
             \textfont1=\scriptfont1
             \scriptfont1=\scriptscriptfont1
             \textfont2=\scriptfont2
             \scriptfont2=\scriptscriptfont2
             \textfont3=\scriptfont3
             \scriptfont3=\scriptscriptfont3
           \renewcommand{\arraystretch}{0.7}
           \begin{array}{l}}{\end{array}}

\newenvironment{scarray}{
             \textfont0=\scriptfont0
             \scriptfont0=\scriptscriptfont0
             \textfont1=\scriptfont1
             \scriptfont1=\scriptscriptfont1
             \textfont2=\scriptfont2
             \scriptfont2=\scriptscriptfont2
             \textfont3=\scriptfont3
             \scriptfont3=\scriptscriptfont3
           \renewcommand{\arraystretch}{0.7}
           \begin{array}{c}}{\end{array}}

\clearpage

One of the pillars of functional analysis is the
uniform boundedness theorem:

\bigskip

{\sc Uniform boundedness theorem.}
Let $\scrf$ be a family of bounded linear operators
from a Banach space $X$ to a normed linear space $Y$.
If $\scrf$ is pointwise bounded
(i.e., $\sup_{T \in \scrf} \|Tx\| < \infty$ for all $x \in X$),
then $\scrf$ is norm-bounded
(i.e., $\sup_{T \in \scrf} \|T\| < \infty$).

\bigskip

The standard textbook proof
(e.g., \cite[p.~81]{Reed_72}),
which goes back to Stefan Banach, Hugo Steinhaus,
and Stanis\l{}aw Saks in 1927 \cite{Banach_27},
employs the Baire category theorem or some variant thereof.\footnote{
   See \cite[p.~319, note~67]{Birkhoff_84} concerning credit to Saks.
}
This proof is very simple,
but its reliance on the Baire category theorem makes it
not completely elementary.

By contrast, the original proofs
given by Hans Hahn \cite{Hahn_22} and Stefan Banach \cite{Banach_22}
in 1922 were quite different:
they began from the assumption that $\sup_{T \in \scrf} \|T\| = \infty$
and used a ``gliding hump'' (also called ``sliding hump'') technique
to construct a sequence $(T_n)$ in $\scrf$ and a point $x \in X$
such that $\lim_{n \to\infty} \|T_n x \| = \infty$.\footnote{
   Hahn's proof is discussed in at least two modern textbooks:
   see \cite[Exercise~1.76, p.~49]{Megginson_98} and
   \cite[Exercise~3.15, pp.~71--72]{MacCluer_09}.
}
Variants of this proof were later given by
T.~H.~Hildebrandt \cite{Hildebrandt_23}
and Felix Hausdorff \cite{Hausdorff_32,Hennefeld_80}.\footnote{
   See also \cite[pp.~63--64]{Riesz_55} and \cite[pp.~74--75]{Weidmann_80}
   for an elementary proof that is closely related to the standard
   ``nested ball'' proof of the Baire category theorem;
   and see \cite[Problem 27, pp.~14--15 and 184]{Halmos_67}
   and \cite{Holland_69} for elementary proofs
   in the special case of linear functionals on a Hilbert space.
}
These proofs are elementary, but the details are a bit fiddly.

Here is a {\em really}\/ simple proof
along similar lines:

\bigskip

{\bf Lemma.}   Let $T$ be a bounded linear operator
from a normed linear space $X$ to a normed linear space $Y$.
Then for any $x \in X$ and $r > 0$, we have
\be
   \sup\limits_{x' \in B(x,r)} \| Tx' \| \;\ge\; \|T\| r  \;,
\ee
where $B(x,r) = \{x' \in X \colon\: \|x'-x\| < r \}$.

\smallskip

\proof
For $\xi \in X$ we have
\be
   \max\bigl\{ \| T(x+\xi) \| ,\, \| T(x-\xi) \| \bigr\}
   \;\,\ge\;\,
   \smhalf \bigl[ \| T(x+\xi) \| + \| T(x-\xi) \| \bigr]
   \;\,\ge\;\,
   \| T \xi \|  \;,
   \quad
\ee
where the second $\ge$ uses the triangle inequality
in the form $\| \alpha-\beta \| \le \|\alpha\| + \|\beta\|$.
Now take the supremum over $\xi \in B(0,r)$.
\qed

\vspace*{-4mm}

\proofof{the uniform boundedness theorem}
Suppose that $\sup_{T \in \scrf} \|T\| = \infty$,
and choose $(T_n)_{n=1}^\infty$ in $\scrf$ such that $\|T_n\| \ge 4^n$.
Then set $x_0 = 0$, and for $n \ge 1$ use the lemma to
choose inductively $x_n \in X$
such that $\| x_n - x_{n-1} \| \le 3^{-n}$
and $\| T_n x_n \| \ge \smfrac{2}{3} 3^{-n} \| T_n \|$.
The sequence $(x_n)$ is Cauchy, hence convergent to some $x \in X$;
and it is easy to see that
$\| x-x_n \| \le \smfrac{1}{2} 3^{-n}$ and hence
$\| T_n x \| \ge \smfrac{1}{6} 3^{-n} \| T_n \| \ge \smfrac{1}{6} (4/3)^n
 \to \infty$.
\qed

\medskip

{\bf Remarks.}
1.  As just seen, this proof is most conveniently expressed in terms of a
{\em sequence}\/ $(x_n)$ that converges to $x$.
This contrasts with the earlier ``gliding hump'' proofs,
which used a {\em series}\/ that sums to $x$.
Of course, sequences and series are equivalent,
so each proof can be expressed in either language;
it is a question of taste which formulation one finds simpler.

2.  This proof is extremely wasteful from a quantitative point of view.
A quantitatively sharp version of the uniform boundedness theorem
follows from Ball's ``plank theorem'' \cite{Ball_91}:
namely, if $\sum_{n=1}^\infty \|T_n\|^{-1} < \infty$,
then there exists $x \in X$ such that
$\lim_{n\to\infty} \| T_n x \| = \infty$
(see also \cite{Muller_09}).

3. A similar (but slightly more complicated) elementary proof
of the uniform boundedness theorem can be found in \cite[p.~83]{Gohberg_03}.

4.  ``Gliding hump'' proofs continue to be useful in functional analysis:
see \cite{Swartz_96} for a detailed survey.

5.  The standard Baire category method
yields a slightly stronger version of the uniform boundedness theorem
than the one stated here, namely:
if $\sup_{T \in \scrf} \|Tx\| < \infty$
for a {\em nonmeager}\/ (i.e., second category) set of $x \in X$,
then $\scrf$ is norm-bounded.

6.  The uniform boundedness theorem has generalizations
to suitable classes of non-normable
and even non-metrizable
topological vector spaces
(see, e.g., \cite[pp.~82--87]{Schaefer_99}).
I~leave it to others to determine whether any ideas from this proof
can be carried over to these more general settings.

7.  More information on the history of the uniform boundedness theorem
can be found in
\cite[pp.~302, 319n67]{Birkhoff_84},
\cite[pp.~138--142]{Dieudonne_81}, and
\cite[pp.~21--22, 40--43]{Pietsch_07}.

\bigskip
\bigskip

\noindent
{\bf Acknowledgments.}
I wish to thank Keith Ball for reminding me of Hahn's ``gliding hump'' proof;
it~was my attempt to fill in the details of Keith's sketch
that led to the proof reported here.
Keith informs me that versions of this proof have been independently devised
by at least four or five people,
though to his knowledge none of them have bothered to publish it.
I also wish to thank David Edmunds and Bob Megginson
for helpful correspondence, and three anonymous referees
for valuable suggestions concerning the exposition.
Finally, I wish to thank J\"urgen Voigt and Markus Haase
for drawing my attention to the elementary proofs in
\cite{Gohberg_03,Riesz_55,Weidmann_80}.

This research was supported in part by
U.S.\ National Science Foundation grant PHY--0424082.

\bigskip
\bigskip



\end{document}